\documentclass{article}

\usepackage{amsmath,amssymb,amsthm,graphicx,proof,color}
\usepackage[UKenglish]{babel}

\usepackage{MnSymbol}
\usepackage{hyperref}
\usepackage[all]{xy}

\DeclareSymbolFont{symbolsC}{U}{txsyc}{m}{n}
\DeclareMathSymbol{\strictif}{\mathrel}{symbolsC}{74}

\newcommand{\bis}{\mathrel{\mathchoice%
{\raisebox{.3ex}{$\,
  \underline{\makebox[.7em]{$\leftrightarrow$}}\,$}}%
{\raisebox{.3ex}{$\,
  \underline{\makebox[.7em]{$\leftrightarrow$}}\,$}}%
{\raisebox{.2ex}{$\,
  \underline{\makebox[.5em]{\scriptsize$\leftrightarrow$}}\,$}}%
{\raisebox{.2ex}{$\,
  \underline{\makebox[.5em]{\scriptsize$\leftrightarrow$}}\,$}}}}

\hypersetup{
   colorlinks,%
   citecolor=blue,%
   filecolor=black,%
   linkcolor=red,%
   urlcolor=black
 }

\usepackage[all]{xy}
\usepackage{tikz}
\usetikzlibrary{trees}
\usepgflibrary{arrows}
\usetikzlibrary{positioning,shapes}
\usetikzlibrary{patterns}

\newcommand{\trianglecdot}[1][0.5]{%
\begin{tikzpicture}
\node(t)[black,regular polygon,regular polygon sides=3,rotate=0,draw,scale=#1] at
 (0,0){};
\node at (t.center)[circle,fill=black,scale=0.5*#1]{};
\end{tikzpicture}}

\newcommand{\downtrianglecdot}[1][0.5]{%
\begin{tikzpicture}
\node(t)[black,regular polygon,regular polygon sides=3,rotate=180,draw,scale=#1] at
 (0,0){};
\node at (t.center)[circle,fill=black,scale=0.5*#1]{};
\end{tikzpicture}}

\newcommand{\M}{\mathcal{M}}
\newcommand{\N}{\mathcal{N}}

\newcommand{\BP}{\ensuremath{\textbf{P}}}

\newcommand{\lr}[1]{\langle #1 \rangle}
\newcommand{\lra}{\leftrightarrow}

\theoremstyle{definition}

\renewcommand{\phi}{\varphi}

\newtheorem{theorem}{Theorem}
\newtheorem{lemma}[theorem]{Lemma}
\newtheorem{definition}[theorem]{Definition}
\newtheorem{proposition}[theorem]{Proposition}

\newtheorem{corollary}[theorem]{Corollary}
\newtheorem{remark}[theorem]{Remark}
\newtheorem{fact}[theorem]{Fact}

\newcommand{\weg}[1]{}

\title{A logic of strong contingency}

\author{Jie Fan\\
\small Institute of Philosophy, Chinese Academy of Sciences;\\
\small School of Humanities, University of Chinese Academy of Sciences\\
\small \texttt{jiefan@ucas.ac.cn}}
\date{}

\begin{document}
\maketitle

\begin{abstract}

In this paper, inspired by the notions of contingency and strong possibility, we propose a notion which may be called `strong contingency'. This notion can be seen as a mixture of contingency and strong possibility. We investigate the logic with strong contingency as a sole primitive modality, including its bisimulation notion, expressivity, frame definability, and axiomatizations. It turns out that this logic is expressively incomparable with the standard modal logic and the logic of contingency even over the class of $\mathcal{S}5$-models, and less expressive than the logic of strong possibility and weak necessity over the class of $\mathcal{S}5$-models. Besides, as the introduction of the inequality into the semantics makes the completeness proof much harder, instead of showing the completeness via a canonical model construction, we reduce the soundness and strong completeness of the involved proof systems to those of the existing proof systems in the literature.

    
\end{abstract}


\section{Introduction}

As is known, standard modal logic concerns the notions of necessity and possibility. Semantically, a formula is necessary, if it is true at all accessible (possible) worlds; a formula is possible, if it is true at some accessible worlds. For standard modal logic, we refer to any textbook on modal logic (e.g., ~\cite {blackburnetal:2001}). Following the tradition of standard modal logic, we use $\Box$ and $\Diamond$ for the operators of necessity and possibility, respectively.

On the one hand, the notion of possibility concerns one-sided possibility, namely, possibly true. In contrast, if considering two-sided possibility, that is, possibly true and possibly false, then we obtain the notion of contingency~\cite{Borgan67,Malink:2016}. Semantically, a formula is contingent, if it is true at some accessible worlds and also false at some accessible worlds; otherwise, it is non-contingent, that is, if all accessible worlds agree on the truth value of the formula. There have been a bunch of logical studies on contingency and non-contingency, see e.g.~\cite{MR66,Cresswell88,Humberstone95,DBLP:journals/ndjfl/Kuhn95,DBLP:journals/ndjfl/Zolin99,hoeketal:2004,Fanetal:2014,Fanetal:2015,Fan:2019symmetric}. Following the tradition of (non-)contingency logic, we use $\Delta$ and $\nabla$ for the operators of non-contingency and contingency, respectively.

On the other hand, the relational semantics of necessity and possibility refers not only to non-actual worlds (that is, other worlds), but to the actual world, namely the evaluated world. The notion of non-actual worlds is related to the actualism vs. possibilism dispute and counterfactual reasoning and imagination, see~\cite{Fan:2023} and the references therein. If considering the non-actual worlds only, we obtain the notions of weak necessity and strong possibility. Semantically, a formula is weakly necessary, if it is true at all accessible non-actual worlds; a formula is strongly possible, if it is true at some accessible non-actual worlds. For the logical studies of strong possibility and weak necessity, we refer to~\cite{Chen:2020,Fan:2023}. Following the tradition of the logic of strong possibility and weak necessity in~\cite{Chen:2020}, alias the logic of non-actual possible worlds in~\cite{Fan:2023}, we use $\boxdot$ and $\diamonddot$ for the operators of weak necessity and strong possibility, respectively.

Inspired by the notions of contingency and strong possibility, we propose a notion which may be called `strong contingency'. Semantically, a formula is {\em strongly contingent}, if it is true at some accessible other worlds and also false at some accessible other worlds; otherwise, it is {\em weakly non-contingent}. The notion of strong contingency can be seen as a mixture of contingency and strong possibility: as in the case of contingency, it considers two-sided possibility; as in the case of strong possibility, it only refers to non-actual worlds. In what follows, we will use $\downtrianglecdot$ and $\trianglecdot$ for the operators of strong contingency and weak non-contingency, respectively.



The notion of weak necessity is related to the notion of elsewhere or elsewhen, in the sense that the latter is a special instance of the former when the accessibility relation is the universal relation (see Fn.~\ref{note.elsewhere}), in which case $\boxdot \phi$ means ``$\phi$ is true everywhere else (or true everywhen else)''~\cite{Segerberg:1976,Wright:1979,Segerberg:1980,deRijke:1992,Humberstone:2007}. Similarly, if the accessibility relation is the universal relation, $\trianglecdot \phi$ means ``any two worlds other than this one agree on the truth-value of $\phi$.'' This hints that there are some relationships between the notions of weak necessity and of weak non-contingency. Indeed, as we will see, the relationship between the two notions is analogous to that between the notions of necessity and of non-contingency.


This allows us to ask whether, when two worlds (in different models on the same frame, perhaps) verify the same $\trianglecdot$ formulas, they must verify the same $\Box$-formulas. This question would be trivial for $\Delta$-formulas and $\Box$-formulas in any modal logic with the T-axiom, since in that case we have an explicit definition of $\Box$ in terms of $\Delta$ (and conjunction), namely, $\Box\phi=_{df}\Delta\phi\land\phi$, as in the Routley-Montgomery noncontingency-based modal logics of the 1960s~\cite{MR66,MR68}. But here, we don't have the automatic equivalence of $\Box \phi$ (or $\boxdot\phi$) with $\trianglecdot\phi \land \phi$. As a matter of fact, as we will show below, neither $\Box$ nor $\boxdot$ is definable with $\trianglecdot$ even over the class of $\mathcal{S}5$-models. This is a significant difference between $\trianglecdot$ and $\Delta$.

In this paper, we investigate the logic of strong contingency. The remainder is structured as follows. Section~\ref{sec.syntaxsemantics} proposes the languages of the logic of strong contingency $\mathcal{L}(\trianglecdot)$ and related logics and their semantics. After that, we define the suitable notion of bisimulation for $\mathcal{L}(\trianglecdot)$ in Section~\ref{sec.bisimulation}, and then investigate the expressivity (Section~\ref{sec.expressivity}) and frame definability (Section~\ref{sec.framedefinability}) of $\mathcal{L}(\trianglecdot)$. To completely axiomatize $\mathcal{L}(\trianglecdot)$ over various frame classes via a reduction method (Section~\ref{sec.axiomatizations-trianglecdot}), we introduce a notion of irreflexive reduction in Section~\ref{sec.irreflexivereduction}, and review the axiomatizations of $\mathcal{L}(\Delta)$ over various frame classes in Section~\ref{sec.axiomatizations-delta}. Finally, we conclude with some future work in Section~\ref{sec.conclusion}.


\section{Syntax and Semantics}\label{sec.syntaxsemantics}

Throughout the paper, we let $\BP$ be a nonempty set of propositional variables. In what follows, we introduce the language of the logics involved in this paper, although we will mainly focus on the logic of strong contingency. 
\begin{definition} 
The language $\mathcal{L}$ is defined recursively as follows:
\[
\begin{array}{lll}
   \phi & := & p\mid \neg\phi\mid (\phi\land\phi)\mid \Box\phi \mid \Delta\phi \mid \boxdot\phi\mid \trianglecdot\phi,\\ 
\end{array}
\]
where $p\in\BP$. With the construct $\Box\phi$ only, we obtain the language $\mathcal{L}(\Box)$ of standard modal logic (see e.g.~\cite{blackburnetal:2001}); with the construct $\Delta\phi$ only, we obtain the language $\mathcal{L}(\Delta)$ of contingency logic (see e.g.~\cite{MR66,Humberstone95}); with the construct $\boxdot\phi$ only, we obtain the language $\mathcal{L}(\boxdot)$ of the logic of non-actual possible worlds in~\cite{Fan:2023}, called ``the logic of strong possibility and weak necessity'' in~\cite{Chen:2020}; with the construct $\trianglecdot\phi$ only, we obtain {\em the language $\mathcal{L}(\trianglecdot)$ of the logic of strong contingency}.  
\end{definition}

Intuitively, $\Box\phi$ is read ``it is necessary that $\phi$'', $\Delta\phi$ is read ``it is non-contingent that $\phi$'', $\boxdot\phi$ is read ``it is weakly necessary that $\phi$'', and $\trianglecdot\phi$ is read ``it is {\em weakly non-contingent} that $\phi$''. Other connectives are defined as usual. In particular, we use $\Diamond\phi$ (read ``it is possible that $\phi$''), $\nabla\phi$ (read ``it is contingent that $\phi$''), $\diamonddot\phi$ (read ``it is strongly possible that $\phi$''), and $\downtrianglecdot\phi$ (read ``it is strongly contingent that $\phi$'') as the abbreviations of $\neg{\Box\neg\phi}$, $\neg\Delta\phi$, $\neg{\boxdot\neg\phi}$, and $\neg\trianglecdot\phi$, respectively.

The language $\mathcal{L}$ is interpreted on (Kripke) models. A {\em model} is a triple $\M=\lr{W,R,V}$, where $W$ is a nonempty set of possible worlds, $R$ is a binary relation over $W$, called ``accessibility relation'', and $V$ is a valuation. A {\em frame} is a model without a valuation. A {\em pointed model} is a pair consisting of a model and a world within it. We use $\mathcal{S}5$-model to denote reflexive and Euclidean model. We use $\mathbb{K}$ and $\mathbb{S}5$ to denote the class of all models and the class of reflexive and Euclidean models, respectively. Moreover, we use $K$, $T$, $B$, and $TB$, respectively, to denote the class of all frames, the class of reflexive frames, the class of symmetric frames, and the class of reflexive and symmetric frames.

\begin{definition} Given a pointed model $(\M,w)$, the semantics of $\mathcal{L}$ is recursively defined as follows.\footnote{Note that in the clause $\boxdot\phi$, if $R$  is universal (that is, $wRu$ for all $w,u\in W$), then we obtain the notion of elsewhere or elsewhen in the literature, see ~\cite{Segerberg:1976,Wright:1979,Segerberg:1980,deRijke:1992,Humberstone:2007}.\label{note.elsewhere}}
\[
\begin{array}{|lll|}
\hline
\M,w\vDash p & \iff & w\in V(p)\\
\M,w\vDash\neg\phi & \iff & \M,w\nvDash\phi \\
\M,w\vDash\phi\land\psi & \iff & \M,w\vDash\phi\text{ and }\M,w\vDash\psi\\
\M,w\vDash \Box\phi & \iff & \text{for every }u\in W,\text{ if }wRu,\text{ then }\M,u\vDash\phi\\
\M,w\vDash \Delta\phi & \iff & \text{for every }u,v\in W,\text{ if }wRu\text{ and }wRv,\\
    & & \text{then }(\M,u\vDash\phi\iff \M,v\vDash\phi)\\
\M,w\vDash\boxdot\phi & \iff & \text{for every }u\in W,\text{ if }wRu\text{ and }w\neq u, \text{ then }\M,u\vDash\phi\\
  \M,w\vDash\trianglecdot\phi & \iff  & \text{for every }u,v\in W,\text{ if }wRu\text{ and }w\neq u\text{ and }wRv\text{ and }w\neq v,\\
    & & \text{then }(\M,u\vDash\phi\iff \M,v\vDash\phi)\\
\hline
\end{array}
\]

We say that $\phi$ is {\em true} at $w$ in $\M$, if $\M,w\vDash\phi$; $\phi$ is {\em valid on a model $\M$}, denoted $\M\vDash\phi$, if for all $w$ in $\M$, $\M,w\vDash\phi$; $\phi$ is {\em valid over a class of models $\mathbb{M}$}, denoted $\mathbb{M}\vDash\phi$, if for all $\M$ in $\mathbb{M}$, $\M\vDash\phi$; $\phi$ is {\em valid on a frame $\mathcal{F}$}, denoted $\mathcal{F}\vDash\phi$, if for all models $\M$ based on $\mathcal{F}$, $\M\vDash\phi$; $\phi$ is {\em valid over a class of frames $F$}, denoted $F\vDash\phi$, if for all frames $\mathcal{F}$ in $F$, $\mathcal{F}\vDash\phi$, and when $F$ is the class of all frames, we also say that $\phi$ is {\em valid}.
\end{definition}

One may easily compute the semantics of other operators.
\[
\begin{array}{lll}
  \M,w\vDash\Diamond\phi & \iff & \text{for some }u\in W,~wRu\text{ and }\M,u\vDash\phi \\
  \M,w\vDash \nabla\phi & \iff & \text{for some }u,v\in W,~wRu\text{ and }wRv\text{ and}\\
  && \M,u\vDash\phi\text{ and }\M,v\nvDash\phi\\
  \M,w\vDash\diamonddot\phi & \iff & \text{for some }u\in W,~wRu\text{ and }w\neq u\text{ and }\M,u\vDash\phi\\
  \M,w\vDash \downtrianglecdot\phi & \iff & \text{for some }u\in W,~wRu\text{ and }w\neq u\text{ and }\M,u\vDash\phi\text{ and}\\
  && \text{for some }v\in W,~wRv\text{ and }w\neq v\text{ and }\M,v\nvDash\phi\\
\end{array}
\]

One can check that $\vDash \diamonddot\phi\to \Diamond\phi$ but $\nvDash\Diamond\phi\to\diamonddot\phi$, and $\vDash\Box\phi\to \boxdot\phi$ but $\nvDash \boxdot\phi\to\Box\phi$. So $\diamonddot$ is deductively stronger than $\Diamond$, and $\boxdot$ is deductively weaker than $\Box$. That is why the terms `strong possibility' and `weak necessity' are used for $\diamonddot$ and $\boxdot$, respectively. Similarly, as $\vDash \downtrianglecdot\phi\to \nabla\phi$ but $\nvDash \nabla\phi\to\downtrianglecdot\phi$, and $\vDash \Delta\phi\to\trianglecdot\phi$ but $\nvDash\trianglecdot\phi\to\Delta\phi$, we use the terms `strong contingency' and `weak non-contingency' for $\downtrianglecdot$ and $\trianglecdot$, respectively.

We can also adopt the following equivalent semantics of $\Delta$ and $\trianglecdot$:
\[
\begin{array}{lll}
\M,w\vDash\Delta\phi & \iff  & \text{{\em either} for every }u\in W,\text{ if }wRu,\text{ then }\M,u\vDash\phi \\
    & & \text{{\em or} for every }u\in W,\text{ if }wRu,\text{ then }\M,u\nvDash\phi.\\
   \M,w\vDash\trianglecdot\phi & \iff  & \text{{\em either} for every }u\in W,\text{ if }wRu\text{ and }w\neq u,\text{ then }\M,u\vDash\phi \\
    & & \text{{\em or} for every }u\in W,\text{ if }wRu\text{ and }w\neq u,\text{ then }\M,u\nvDash\phi.\\
\end{array}
\]

One may easily see that although the semantics of $\trianglecdot$ and $\Delta$ are different in general, they are the same on the irreflexive models, that is, the models $\lr{W,R,V}$ in which the following condition is satisfied:
\begin{center}
    for all $x,y\in W$, if $xRy$, then $x\neq y$.
\end{center}

Recall that $\Delta$ is definable with $\Box$ in general, as $\vDash\Delta\phi\lra \Box\phi\vee\Box\neg\phi$, but the converse does not hold, even though $\Box$ is indeed definable with $\Delta$ on reflexive models, as $\Box\phi\lra \Delta\phi\land\phi$ is valid over the class of reflexive models. Despite this, $\Box$ is almost definable in terms of $\Delta$, which means that $\Box$ is definable with $\Delta$ given some precondition; more specifically, $\neg\Delta\psi\to (\Box\phi\lra \Delta\phi\land \Delta(\psi\to\phi))$ is a validity. All this also applies to the case for $\trianglecdot$ and $\boxdot$. Indeed, by semantics, one should easily verify that $\vDash \trianglecdot\phi\lra \boxdot\phi\vee\boxdot\neg\phi$, thus $\trianglecdot$ is definable with $\boxdot$ in general. Also, $\boxdot$ is not definable with $\trianglecdot$ in general --- as a matter of fact, as we shall show, the undefinability result even goes to $\mathcal{S}5$-models. In spite of this, $\boxdot$ is almost definable in terms of $\trianglecdot$, as shown below.

\begin{proposition}\label{prop.almost-definability}
    $\vDash \neg\trianglecdot\psi\to (\boxdot\phi\lra \trianglecdot\phi\land\trianglecdot(\psi\to\phi))$.
\end{proposition}

\begin{proof}
    Let $\M=\lr{W,R,V}$ be a model and $w\in W$. Suppose that $\M,w\vDash\neg\trianglecdot\psi$, to show that $\M,w\vDash \boxdot\phi\lra \trianglecdot\phi\land\trianglecdot(\psi\to\phi)$.

    Assume that $\M,w\vDash \boxdot\phi$. From the truth of $\boxdot\phi$ at $w$ directly follows not only the truth of $\trianglecdot\phi$ but also the truth of $\trianglecdot(\psi\to\phi)$ at $w$ (as all non-actual accessible worlds of $w$ satisfy $\phi$, they also satisfy $\psi\to \phi$, which is weaker). 

    Conversely, assume, for reductio, that $\M,w\vDash \trianglecdot\phi\land \trianglecdot(\psi\to\phi)$, but $\M,w\nvDash\boxdot\phi$. Then there exists $u\in W$ such that $wRu$ and $w\neq u$ and $\M,u\nvDash\phi$. By supposition, there are $v_1,v_2$ such that $wRv_1$ and $wRv_2$ and $w\neq v_1$ and $w\neq v_2$ and $\M,v_1\vDash \psi$ and $\M,v_2\nvDash\psi$, and thus $\M,v_2\vDash \psi\to\phi$. Then using $\M,w\vDash\trianglecdot(\psi\to\phi)$, we can obtain that $\M,v_1\vDash\psi\to\phi$, and thus $\M,v_1\vDash\phi$. Then $\M,w\nvDash\trianglecdot\phi$: a contradiction.
\end{proof}

\section{Bisimulation}\label{sec.bisimulation}

Recall that the bisimulation notions for $\mathcal{L}(\Box)$, $\mathcal{L}(\Delta)$, and $\mathcal{L}(\boxdot)$ have been proposed in the literature, see~\cite{blackburnetal:2001}, \cite{Fanetal:2014}, and~\cite{Chen:2020}, respectively. As we shall see in Section~\ref{sec.expressivity}, however, none of these notions is suitable for $\mathcal{L}(\trianglecdot)$. In this section, we provide a suitable notion of bisimulation for $\mathcal{L}(\trianglecdot)$.

\begin{definition}
    Let $\M=\lr{W,R,V}$ be a model. A nonempty binary relation $Z$ is a {\em $\trianglecdot$-bisimulation} on $\M$, if $wZw'$ implies the following holds:
    \begin{itemize}
        \item[(Var)] for all $p\in\BP$, $w\in V(p)$ iff $w'\in V(p)$;
        \item[($\trianglecdot$-Zig)] if $wRx$ and $w\neq x$, then if there are $u,v$ such that $wRu$ and $wRv$ and $w\neq u$ and $w\neq v$ and $(u,v)\notin Z$, then there exists $x'$ such that $w'Rx'$ and $w'\neq x'$ and $xZx'$;
        \item[($\trianglecdot$-Zag)] if $w'Rx'$ and $w'\neq x'$, then if there are $u',v'$ such that $w'Ru'$ and $w'Rv'$ and $w'\neq u'$ and $w'\neq v'$ and $(u',v')\notin Z$, then there exists $x$ such that $wRx$ and $w\neq x$ and $xZx'$.
    \end{itemize}
    We say that $(\M,w)$ and $(\M',w')$ are {\em $\trianglecdot$-bisimilar}, denoted $(\M,w)\bis_{\trianglecdot}(\M',w')$, if there is a $\trianglecdot$-bisimulation $Z$ on the disjoint union of $\M$ and $\M'$ such that $wZw'$.
\end{definition}

This definition is inspired by the almost definability schema of $\boxdot$ in terms of $\trianglecdot$ (Prop.~\ref{prop.almost-definability}). Recall that the schema says that $\boxdot$ is definable with $\trianglecdot$, given a condition $\neg\trianglecdot\psi$ for some $\psi$. This condition corresponds to a precondition that the current world can see two non-$\mathcal{L}(\trianglecdot)$-equivalent successors distinct from it. Thus we can revise the $\boxdot$-Zig and $\boxdot$-Zag conditions in the definition of $\diamonddot$-bisimulation\footnote{The bisimulation of $\mathcal{L}(\boxdot)$, called `$\diamonddot$-bisimulation' in~\cite[Def.~1]{Chen:2020}, is defined as follows (we revise the definition a little bit to make the comparison clearer).

Let $\M=\lr{W,R,V}$ and $\M'=\lr{W',R',V'}$ be models. A binary relation $Z\subseteq W\times W'$ is said to be a {\em $\diamonddot$-bisimulation} between $\M$ and $\M'$, denoted $Z:\M\bis_{\diamonddot}\M'$, if $wZw'$ implies the following conditions holds:
  \begin{itemize}
    \item[(Var)] for all $p\in\BP$, $w\in V(p)$ iff $w'\in V'(p)$;
    \item[($\boxdot$-Zig)] if $wRx$ and $w\neq x$, then there exists $x'$ such that $w'R'x'$ and $w'\neq x'$ and $xZx'$;
    \item[($\boxdot$-Zag)] if $w'R'x'$ and $w'\neq x'$, then there exists $x$ such that $wRx$ and $w\neq x$ and $xZx'$.
  \end{itemize} 

We say that $(\M,w)$ and $(\M',w')$ are {\em $\diamonddot$-bisimilar}, denoted $(\M,w)\bis_{\diamonddot}(\M',w')$, if there is a $\diamonddot$-bisimulation $Z$ between $\M$ and $\M'$ such that $(w,w')\in Z$.
} in~\cite{Chen:2020} to accommodate the semantic condition, so that we obtain the $\trianglecdot$-Zig and $\trianglecdot$-Zag conditions. Note that for technical convenience, we define $\trianglecdot$-bisimulation within a single model, since the new Zig and Zag conditions require a precondition about ‘sibling’ worlds that are distinct from the current world, i.e. the structural counterpart of non-$\trianglecdot$-equivalence. Based on $\trianglecdot$-bisimulation, we define $\trianglecdot$-bisimilarity between different models. This is similar to the case of $\Delta$-bisimulation and $\Delta$-bisimilarity~\cite[Def.~3.3]{Fanetal:2014}.

We claim without proof that $\trianglecdot$-bisimilarity is an equivalence relation. The following result indicates that $\mathcal{L}(\trianglecdot)$-formulas are invariant under $\trianglecdot$-bisimilarity; in other words, $\trianglecdot$-bisimilarity implies $\trianglecdot$-equivalence.


\begin{proposition}\label{prop.invariance}
    Let $\M=\lr{W,R,V}$ and $\M'=\lr{W',R',V'}$, and $w\in W$ and $w'\in W'$. If $(\M,w)\bis_{\trianglecdot}(\M',w')$, then $(\M,w)\equiv_{\trianglecdot}(\M',w')$.
\end{proposition}

\begin{proof}
    Suppose that $(\M,w)\bis_{\trianglecdot}(\M',w')$. Then there exists a $\trianglecdot$-bisimulation $Z$ on the disjoint union of $\M$ and $\M'$ such that $wZw'$. We proceed by induction on $\phi$, where the nontrivial case is $\trianglecdot\phi$.

    Assume that $\M,w\nvDash\trianglecdot\phi$. Then there are $u,v$ such that $wRu$ and $wRv$ and $w\neq u$ and $w\neq v$ and $\M,u\vDash \phi$ and $\M,v\nvDash\phi$. By induction hypothesis, $(u,v)\notin Z$. Then by using ($\trianglecdot$-Zig) and $wRu$ and $w\neq u$, we infer that there exists $u'$ such that $w'R'u'$ and $w'\neq u'$ and $uZu'$, and thus $\M',u'\vDash\phi$. Similarly, from $wRv$ and $w\neq v$, we can derive that there exists $v'$ such that $w'R'v'$ and $w'\neq v'$ and $vZv'$, and thus $\M',v'\nvDash\phi$. Hence, $\M',w'\nvDash\trianglecdot\phi$. The converse is analogous.
\end{proof}

Conversely, on image-finite models, $\trianglecdot$-equivalence implies $\trianglecdot$-bisimilarity. A model $\M=\lr{W,R,V}$ is said to be {\em image-finite}, if for all $w\in W$, $R(w)=\{u\in W\mid wRu\}$ is finite.

\begin{theorem}[Hennessy-Milner Theorem for $\mathcal{L}(\trianglecdot)$]
    Let $\M=\lr{W,R,V}$ and $\M'=\lr{W',R',V'}$ be image-finite models, $w\in W$ and $w'\in W'$. Then 
    $$(\M,w)\equiv_{\trianglecdot}(\M',w')\text{ iff }(\M,w)\bis_{\trianglecdot}(\M',w').$$
\end{theorem}

\begin{proof}
    By Prop.~\ref{prop.invariance}, it remains only to show `only if'. For this, we need to show that $\equiv_{\trianglecdot}$ is a $\trianglecdot$-bisimulation on the disjoint union of $\M$ and $\M'$. The condition (Var) is immediate. In what follows, we only show that $\equiv_{\trianglecdot}$ satisfies ($\trianglecdot$-Zig), as the proof of ($\trianglecdot$-Zag) is similar.

    Suppose that $w\equiv_{\trianglecdot} w'$, $wRu$ and $wRv$ and $w\neq u$ and $w\neq v$ and $u\not\equiv_{\trianglecdot} v$ and $wRx$ and $w\neq x$, to find an $x'$ such that $w'R'x'$ and $w'\neq x'$ and $x\equiv_{\trianglecdot} x'$. Let $S'=\{u'\mid w'R'u' \text{ and }w'\neq u'\}$. By $u\not\equiv_{\trianglecdot} v$, there exists $\phi\in\mathcal
    {L}(\trianglecdot)$ such that $\M,u\vDash\phi$ and $\M,v\nvDash\phi$. Then $\M,w\nvDash\trianglecdot\phi$. Therefore $\M',w'\nvDash\trianglecdot\phi$, so that $S'$ must be nonempty. Also, as $\M'$ is image-finite, $S'$ must be finite, say $S'=\{x_1',\dots,x_n'\}$. If for each $x_i'\in S'$ there exists $\psi_i\in\mathcal{L}(\trianglecdot)$ such that $\M,x\vDash \psi_i$ and $\M',x_i'\nvDash\psi_i$, then $\M,x\vDash\psi_1\land\cdots\land\psi_n$ and $\M',x_i'\nvDash\psi_1\land\cdots\land\psi_n$, thus $\M',x_i'\vDash\psi_1\land\cdots\land\psi_n\to\phi$. Hence $\M',w'\vDash\trianglecdot(\psi_1\land\cdots\land\psi_n)$ and $\M',w'\vDash\trianglecdot(\psi_1\land\cdots\land\psi_n\to\phi)$. Using supposition again, we infer that $\M,w\vDash\trianglecdot(\psi_1\land\cdots\land\psi_n)$ and $\M,w\vDash\trianglecdot(\psi_1\land\cdots\land\psi_n\to\phi)$. From $\M,w\vDash\trianglecdot(\psi_1\land\cdots\land\psi_n)$ and $\M,x\vDash\psi_1\land\cdots\land\psi_n$, $w\neq x$ and $wRu$ and $wRv$ and $w\neq u$ and $w\neq v$, it follows that $\M,u\vDash \psi_1\land\cdots\land\psi_n$ and $\M,v\vDash \psi_1\land\cdots\land\psi_n$. Then $\M,u\vDash\psi_1\land\cdots\land\psi_n\to\phi$ but $\M,v\nvDash\psi_1\land\cdots\land\psi_n\to\phi$, which contradicts $\M,w\vDash\trianglecdot(\psi_1\land\cdots\land\psi_n\to\phi)$. This shows that there exists $\chi'\in S'$ such that $x\equiv_{\trianglecdot}x'$, and therefore there exists $\chi'$ such that $w'R'x'$ and $w'\neq x'$ and $x\equiv_{\trianglecdot}x'$.
\end{proof}

\section{Expressivity and Frame Definability}

In this section, we compare the relative expressivity of the language $\mathcal{L}(\trianglecdot)$ and the languages $\mathcal{L}(\Delta)$, $\mathcal{L}(\Box)$ and $\mathcal{L}(\boxdot)$, and investigate the frame definability of $\mathcal{L}(\trianglecdot)$.

\subsection{Expressivity}\label{sec.expressivity}

We start with the expressivity part. It turns out that $\mathcal{L}(\trianglecdot)$ and $\mathcal{L}(\Delta)$ (and also $\mathcal{L}(\Box)$) are incomparable in expressivity over the class of $\mathcal{S}5$-models, and $\mathcal{L}(\trianglecdot)$ is less expressive than $\mathcal{L}(\boxdot)$ over the class of $\mathcal{S}5$-models. To begin with, we introduce the related notions.
\begin{definition}
  Let two languages $\mathcal{L}_1$ and $\mathcal{L}_2$ be interpreted on the same class of models $\mathbb{M}$.
  \begin{itemize}
    \item $\mathcal{L}_2$ is {\em at least as expressive as }$\mathcal{L}_1$ over $\mathbb{M}$, notation: $\mathcal{L}_2\geq\mathcal{L}_1[\mathbb{M}]$, if for any $\phi\in\mathcal{L}_1$, there exists $\psi\in\mathcal{L}_2$ such that for any $\M$ in $\mathbb{M}$ and any $w$ in $\M$, we have $\M,w\vDash\phi$ iff $\M,w\vDash\psi$.
    \item $\mathcal{L}_2$ is {\em more expressive than }$\mathcal{L}_1$ over $\mathbb{M}$, notation: $\mathcal{L}_2>\mathcal{L}_1[\mathbb{M}]$, if $\mathcal{L}_2\geq \mathcal{L}_1[\mathbb{M}]$ and $\mathcal{L}_1\not\geq\mathcal{L}_2[\mathbb{M}]$. In this case, we also say that $\mathcal{L}_1$ is {\em less expressive than }$\mathcal{L}_2$ over $\mathbb{M}$. 
    \item $\mathcal{L}_1$ and $\mathcal{L}_2$ are {\em incomparable in expressivity} over $\mathbb{M}$, notation: $\mathcal{L}_1\asymp\mathcal{L}_2[\mathbb{M}]$, if $\mathcal{L}_1\not\geq\mathcal{L}_2[\mathbb{M}]$ and $\mathcal{L}_2\not\geq\mathcal{L}_1[\mathbb{M}]$.
  \end{itemize}
\end{definition}

\subsubsection{$\mathcal{L}(\trianglecdot)$ vs. $\mathcal{L}(\Delta)$}

First, $\mathcal{L}(\trianglecdot)$ and $\mathcal{L}(\Delta)$ are incomparable in expressivity over the class of $\mathcal{S}5$-models. This implies that neither $\Delta$ nor $\trianglecdot$ is definable with each other over the class of $\mathcal{S}5$-models.

\begin{proposition}\label{prop.trianglecdot-delta-s5}
    $\mathcal{L}(\trianglecdot)\not\geq \mathcal{L}(\Delta)[\mathbb{S}5]$. That is, over the class of $S5$-models, $\mathcal{L}(\trianglecdot)$ is not at least as expressive as $\mathcal{L}(\Delta)$.
\end{proposition}

\begin{proof}
    Consider the following $S5$-models:
    
    \[
\xymatrix{\M&w:p\ar@(ul,ur)\ar[rr] && u:p\ar[ll]\ar@(ur,ul) && \M'&w':p\ar@(ul,ur)\ar[rr] && u':\neg p\ar@(ur,ul)\ar[ll]}
    \]

    Define $Z=\{(w,w')\}$. Clearly, $(\M,w)\bis_{\trianglecdot}(\M',w')$. By Prop.~\ref{prop.invariance}, $\mathcal{L}(\trianglecdot)$ cannot distinguish between $(\M,w)$ and $(\M',w')$. However, the two pointed models can be distinguished by $\Delta p$, 
 since $\M,w\vDash\Delta p$ and $\M',w'\nvDash\Delta p$. Therefore, $\mathcal{L}(\trianglecdot)\not\geq\mathcal{L}(\Delta)[\mathbb{S}5]$.
\end{proof}

\begin{proposition}\label{prop.delta-trianglecdot}
   $\mathcal{L}(\Delta)\not\geq \mathcal{L}(\trianglecdot)[\mathbb{S}5]$. That is, over the class of $S5$-models, $\mathcal{L}(\Delta)$ is not at least as expressive as $\mathcal{L}(\trianglecdot)$.
\end{proposition}

\begin{proof}
    Consider the following $S5$-models:
    
    \[
    \xymatrix{& & u:p\ar[dl]\ar[dr]\ar@(ul,ur) & \\
    \M& w:p\ar[ur]\ar[rr]\ar@(dl,dr) & & v:\neg p\ar[ll]\ar[ul]\ar@(dr,dl) & \M'& w':p\ar[r]\ar@(ul,ur) & v':\neg p\ar[l]\ar@(ur,ul)\\}
    \]

    \medskip

    Define $Z=\{(w,w'),(u,w'),(v,v')\}$. Obviously, $(\M,w)\bis_{\Delta}(\M',w')$. By~\cite[Prop.~3.5]{Fanetal:2014}, $\mathcal{L}(\Delta)$ cannot distinguish between $(\M,w)$ and $(\M',w')$. However, the two pointed models can be distinguished by $\trianglecdot p$, since $\M,w\nvDash \trianglecdot p$ and $\M',w'\vDash\trianglecdot p$.
\end{proof}

\begin{corollary}
    $\mathcal{L}(\trianglecdot)\asymp \mathcal{L}(\Delta)[\mathbb{S}5]$. That is, over the class of $S5$-models, $\mathcal{L}(\trianglecdot)$ and $\mathcal{L}(\Delta)$ are incomparable in expressivity.
\end{corollary}

\subsubsection{$\mathcal{L}(\trianglecdot)$ vs. $\mathcal{L}(\Box)$}

Second, $\mathcal{L}(\trianglecdot)$ and $\mathcal{L}(\Box)$ are incomparable in expressivity over the class of $\mathcal{S}5$-models. This entails that neither $\Box$ nor $\trianglecdot$ is definable with each other over the class of $\mathcal{S}5$-models.
\begin{proposition}
   $\mathcal{L}(\trianglecdot)\not\geq \mathcal{L}(\Box)[\mathbb{S}5]$. That is, over the class of $S5$-models, $\mathcal{L}(\trianglecdot)$ is not at least as expressive as $\mathcal{L}(\Box)$.
\end{proposition}

\begin{proof}
    Use the models in the proof of  Prop.~\ref{prop.trianglecdot-delta-s5}. Note that $\M,w\vDash\Box p$ and $\M',w'\nvDash\Box p$.
\end{proof}

Thus even over the class of $S5$-models, when two pointed models verify the same $\trianglecdot$-formulas, they do not necessarily verify the same $\Box$-formulas. This entails that $\Box$ is {\em not} definable with $\trianglecdot$ even on the class of $\mathcal{S}5$-models. As mentioned in the introduction, this is unlike the case for $\Delta$, where $\Box$ is definable with $\Delta$ over the class of reflexive models. 

\begin{proposition}
   $\mathcal{L}(\Box)\not\geq \mathcal{L}(\trianglecdot)[\mathbb{S}5]$. That is, over the class of $S5$-models, $\mathcal{L}(\Box)$ is not at least as expressive as $\mathcal{L}(\trianglecdot)$.
\end{proposition}

\begin{proof}
Use the models in the proof of Prop.~\ref{prop.delta-trianglecdot}. Note that $(\M,w)\bis_{\Box}(\M',w')$.
\end{proof}

\begin{corollary}
$\mathcal{L}(\trianglecdot)\asymp \mathcal{L}(\Box)[\mathbb{S}5]$. That is, over the class of $S5$-models, $\mathcal{L}(\trianglecdot)$ and $\mathcal{L}(\Box)$ are incomparable in expressivity.
\end{corollary}

\subsubsection{$\mathcal{L}(\trianglecdot)$ vs. $\mathcal{L}(\boxdot)$}

Third, $\mathcal{L}(\boxdot)$ is at least as expressive as $\mathcal{L}(\trianglecdot)$ over the class of all models, and $\mathcal{L}(\trianglecdot)$ is less expressive than $\mathcal{L}(\boxdot)$ over the class of $\mathcal{S}5$-models. This implies that $\boxdot$ is not definable with $\trianglecdot$ over the class of $\mathcal{S}5$-models, which is contrast to the fact that $\Box$ is definable with $\Delta$ over the class of reflexive models.

\begin{proposition}\label{prop.boxdot>trianglecdot}
    $\mathcal{L}(\boxdot)\geq \mathcal{L}(\trianglecdot)[\mathbb{K}]$. That is, over the class of all models, $\mathcal{L}(\boxdot)$ is at least as expressive as $\mathcal{L}(\trianglecdot)$.
\end{proposition}

\begin{proof}
    This is because $\trianglecdot$ is definable with $\boxdot$, as $\vDash \trianglecdot\phi\lra (\boxdot\phi\vee \boxdot\neg\phi)$.
\end{proof}

\begin{proposition}
    $\mathcal{L}(\boxdot)>\mathcal{L}(\trianglecdot)[\mathbb{S}5]$. That is, over the class of $S5$-models, $\mathcal{L}(\boxdot)$ is more expressive than $\mathcal{L}(\trianglecdot)$.
\end{proposition}

\begin{proof}
    By Prop.~\ref{prop.boxdot>trianglecdot}, it remains to show that $\mathcal{L}(\trianglecdot)\not\geq \mathcal{L}(\boxdot)[\mathbb{S}5]$. Use the models $\M$ and $\M'$ in the proof of Prop.~\ref{prop.trianglecdot-delta-s5}. We have shown that $\mathcal{L}(\trianglecdot)$ cannot distinguish between $(\M,w)$ and $(\M',w')$. However, the two pointed models can be distinguished by $\boxdot p$, because $\M,w\vDash \boxdot p$ and $\M',w'\nvDash\boxdot p$.
\end{proof}

\subsection{Frame definability}\label{sec.framedefinability}

Now we investigate the frame definability of $\mathcal{L}(\trianglecdot)$. It turns out that all five of the most familiar frame properties are undefinable in $\mathcal{L}(\trianglecdot)$. This is like the case in $\mathcal{L}(\Delta)$. 

The following result is immediate by semantics.
\begin{fact}\label{prop.valid-trianglecdotformulas}
    Let $\mathcal{F}=\lr{W,R}$ be a frame. If for all $w\in W$, $|\{u\in W\mid wRu\text{ and }w\neq u\}|<2$, then for all $\phi\in\mathcal{L}(\trianglecdot)$, we have $\mathcal{F}\vDash\trianglecdot\phi$.
\end{fact}

\begin{proposition}
    The frame properties of seriality, reflexivity, transitivity, symmetry, and Euclideanity are not definable in $\mathcal{L}(\trianglecdot)$.
\end{proposition}

\begin{proof}
    Consider the following frames:
    \[\xymatrix{\mathcal{F}:& w\ar[rr]&& u\ar[rr] && v &\mathcal{F}'& w'\ar@(ur,ul)\\
    }\]

    By Fact~\ref{prop.valid-trianglecdotformulas}, we can show $(\ast)$: for all $\phi\in\mathcal{L}(\trianglecdot)$, $\mathcal{F}\vDash\phi$ iff $\mathcal{F}'\vDash\phi$.

    Note that $\mathcal{F}'$ is serial (resp. reflexive, transitive, symmetric, and Euclidean), but $\mathcal{F}$ is not. If seriality (resp. reflexivity, transitivity, symmetry, and Euclideanity) is defined by a set of $\mathcal{L}(\trianglecdot)$-formulas, say $\Gamma$, then, as $\mathcal{F}'$ is serial (resp. reflexive, transitive, symmetric, and Euclidean), we have $\mathcal{F}'\vDash\Gamma$. Due to $(\ast)$, we should also have $\mathcal{F}\vDash\Gamma$, and thus $\mathcal{F}$ is also serial (resp. reflexive, transitive, symmetric, and Euclidean): a contradiction.
\end{proof}


\section{Irreflexive Reduction}\label{sec.irreflexivereduction}

To show the completeness of proof systems over special frames, we need to use a notion of reflexive reduction proposed in~\cite[Def.~2]{Fan:2023}. Intuitively, the irreflexive reduction of a frame is obtained via omitting all reflexive arrows in the original frame. More formally, 
\begin{definition}\label{def.ir}[Irreflexive reduction]
Let $\mathcal{F}=\lr{W,R}$ be a frame. Its {\em irreflexive reduction} is $\mathcal{F}^{-T}=\lr{W,R^{-T}}$, where
$$R^{-T}=R\backslash \{(w,w)\mid w\in W\}.$$
Moreover, $\M^{-T}=\lr{\mathcal{F}^{-T},V}$ is said to be the {\em irreflexive reduction} of $\M=\lr{\mathcal{F},V}$, if $\mathcal{F}^{-T}$ is the irreflexive reduction of $\mathcal{F}$. We say that a class of frames $\mathfrak{C}$ is {\em closed under irreflexivization}, if $\mathcal{F}\in\mathfrak{C}$ implies $\mathcal{F}^{-T}\in\mathfrak{C}$.
\end{definition}

One may easily see that every frame/model has a unique irreflexive reduction.

The satisfiability and frame validity of $\mathcal{L}(\trianglecdot)$-formulas are invariant under the notion of irreflexive reduction.
\begin{proposition}\label{prop.invariant-ir}
    Let $\mathcal{F}^{-T}=\lr{W,R^{-T}}$ be the irreflexive reduction of $\mathcal{F}=\lr{W,R}$, and let $\M=\lr{\mathcal{F},V}$ and $\mathcal{M}^{-T}=\lr{\mathcal{F}^{-T},V}$. Then for all $w\in W$, for all $\phi\in\mathcal{L}(\trianglecdot)$, we have
        $$\M,w\vDash\phi\iff \M^{-T},w\vDash\phi.$$
Consequently, for all $\phi\in\mathcal{L}(\trianglecdot)$, we have
        $\mathcal{F}\vDash\phi\iff \mathcal{F}^{-T}\vDash\phi.$
\end{proposition}

\begin{proof}
    The proof is by induction on $\phi$, where the nontrivial case is $\trianglecdot\phi$.

    Suppose that $\M,w\nvDash\trianglecdot\phi$. Then there are $u,v\in W$ such that $wRu$ and $w\neq u$ and $wRv$ and $w\neq v$ and $\M,u\vDash\phi$ and $\M,v\nvDash\phi$. Clearly, $wR^{-T}u$ and $wR^{-T}v$. By induction hypothesis, we have $\M^{-T},u\vDash\phi$ and $\M^{-T},v\nvDash\phi$. Hence, $\M^{-T},w\nvDash\trianglecdot\phi$.

    Conversely, assume that $\M^{-T},w\nvDash\trianglecdot\phi$. Then there are $u,v\in W$ such that $wR^{-T}u$ and $w\neq u$ and $wR^{-T}v$ and $w\neq v$ and $\M^{-T},u\vDash\phi$ and $\M^{-T},v\nvDash\phi$. Since $R^{-T}\subseteq R$, it follows that $wRu$ and $wRv$. By induction hypothesis, we have $\M,u\vDash\phi$ and $\M,v\nvDash\phi$. Therefore, $\M,w\nvDash\trianglecdot
    \phi$.
\end{proof}

Given a class $\mathfrak{C}$ of frames, define $\mathfrak{C}^{-T}=\{\mathcal{F}^{-T}\mid \mathcal{F}\in \mathfrak{C}\}$, where $\mathcal{F}^{-T}$ is the irreflexive reduction of $\mathcal{F}$. In other words, $\mathfrak{C}^{-T}$ is the set of the irreflexive reduction of each frame in $\mathfrak{C}$.
\begin{lemma}\label{lem.properties-trianglecdot}
    Let $\mathfrak{C}$ and $\mathfrak{C}'$ be two classes of frames. If $\mathfrak{C}^{-T}=\mathfrak{C}'^{-T}$, then for all $\Gamma\cup\{\phi\}\subseteq\mathcal{L}(\trianglecdot)$, we have
    $$\Gamma\vDash_{\mathfrak{C}}\phi\text{ iff }\Gamma\vDash_{\mathfrak{C}'}\phi.$$
\end{lemma}

\begin{proof}
    By Prop.~\ref{prop.invariant-ir}, we can show that $\Gamma\vDash_{\mathfrak{C}}\phi$ iff $\Gamma\vDash_{\mathfrak{C}^{-T}}\phi$, and $\Gamma\vDash_{\mathfrak{C}'}\phi$ iff $\Gamma\vDash_{\mathfrak{C}'^{-T}}\phi$. Since $\mathfrak{C}^{-T}=\mathfrak{C}'^{-T}$, it follows that $\Gamma\vDash_{\mathfrak{C}}\phi\text{ iff }\Gamma\vDash_{\mathfrak{C}'}\phi.$
\end{proof}

Let $X$ be a class of frames. One may check that if $T\subseteq X$, then $X^{-T}=T^{-T}$, and if $TB\subseteq X\subseteq B$, then $TB^{-T}=X^{-T}=B^{-T}$. Then by Lemma~\ref{lem.properties-trianglecdot}, we have the following semantic properties of $\mathcal{L}(\trianglecdot)$.
\weg{\begin{proposition}\label{prop.semanticproperties}
    Let $\Gamma\cup\{\phi\}\subseteq\mathcal{L}(\trianglecdot)$ and $\mathbb{T}\subseteq \mathbb{X}$. Then $\Gamma\vDash_{\mathbb{X}}\phi$ iff $\Gamma\vDash_{\mathbb{T}}\phi$.
    \end{proposition}}
\begin{proposition}\label{prop.semanticproperties}
    Let $\Gamma\cup\{\phi\}\subseteq\mathcal{L}(\trianglecdot)$. Then
    \begin{itemize}
        \item[(1)] Where $T\subseteq X$, $\Gamma\vDash_{X}\phi$ iff $\Gamma\vDash_{T}\phi$;
        \item[(2)] Where $TB\subseteq X\subseteq B$, $\Gamma\vDash_{TB}\phi$ iff $\Gamma\vDash_{X}\phi$ iff $\Gamma\vDash_{B}\phi$.
    \end{itemize}
\end{proposition}

We conclude this section with the following result, which states that $\mathcal{L}(\trianglecdot)$ is insensitive to reflexivity; that is, adding or deleting reflexive arrows does not change the validity/satisfiability of $\mathcal{L}(\trianglecdot)$-formulas in a given frame.
\begin{proposition}
    Let $\M=\lr{W,R,V}$ and $\M'=\lr{W,R',V}$ be models such that $R$ and $R'$ only differ in pairs of reflexive arrows. Then for all $\phi\in\mathcal{L}(\trianglecdot)$, for all $w\in W$, we have
    $$\M,w\vDash \phi\text{ iff }\M',w\vDash\phi.$$
\end{proposition}

\begin{proof}
    Note that $\M$ and $\M'$ must have the same irreflexive reduction. Then use Prop.~\ref{prop.invariant-ir}.
\end{proof}

In what follows, we will completely axiomatize $\mathcal{L}(\trianglecdot)$ over various frame classes. As illustrated by~\cite{Segerberg:1980,deRijke:1992,Humberstone:2007,Chen:2020}, the introduction of the inequality into the semantics makes the completeness proof much harder. Rather than showing the completeness via a canonical model construction, in what follows, we will reduce the soundness and completeness of the minimal system and of the symmetric system of $\mathcal{L}(\trianglecdot)$ to those of the minimal system and of the symmetric system of $\mathcal{L}(\Delta)$. So first we will review the soundness and completeness of the minimal system and the symmetric system of $\mathcal{L}(\Delta)$. 

\section{Axiomatizations of $\mathcal{L}(\Delta)$}\label{sec.axiomatizations-delta}

\begin{definition}\label{def.Delta-systems}
    The proof system ${\bf K^{\Delta}}$ consists of the following axioms and inference rules.
    \[
\begin{array}{ll}
   \text{TAUT}  & \text{all instances of propositional tautologies} \\
   \Delta\text{Con} & \Delta(\chi\to\phi)\land\Delta(\neg\chi\to\phi)\to \Delta\phi\\
   \Delta\text{Dis} & \Delta\phi\to \Delta(\phi\to\psi)\vee\Delta(\neg\phi\to\chi)\\
   \Delta{\lra}  & \Delta\phi\lra \Delta\neg\phi\\
   \text{MP} & \text{From }\phi\text{ and }\phi\to\psi\text{ infer }\psi\\
   \text{NEC}\Delta & \text{From }\phi\text{ infer }\Delta\phi\\
   \text{RE}\Delta & \text{From }\phi\lra \psi\text{ infer }\Delta\phi\lra \Delta\psi\\
\end{array}
\]
The proof system ${\bf KB^{\Delta}}$ is the extension of ${\bf K^{\Delta}}$ with axiom $\Delta\text{B}$: $$\phi\to\Delta((\Delta\phi\land\Delta(\phi\to\psi)\land\neg\Delta\psi)\to\chi).$$
\end{definition}

Notions of derivation and theorems are defined as usual. In particular, we use $\Gamma\vdash_{\bf K^{\Delta}}\phi$ (resp. $\Gamma\vdash_{\bf KB^{\Delta}}\phi$) to denote that $\phi$ is derivable from $\Gamma$ in ${\bf K^{\Delta}}$ (resp. ${\bf KB^{\Delta}}$).

\begin{theorem}\cite{Fanetal:2015}\label{thm.comp-k-delta}
    ${\bf K^{\Delta}}$ is sound and strongly complete with respect to the class of all frames.
\end{theorem}

\begin{theorem}\cite{Fanetal:2014,Fanetal:2015,Fan:2019symmetric}\label{thm.comp-kb-delta}
    ${\bf KB^{\Delta}}$ is sound and strongly complete with respect to the class of symmetric frames.
\end{theorem}

We conclude this section with an important result, which will be used in Lemma~\ref{lem.equivalence}.

\begin{theorem}\label{thm.comp-kdelta-irreflexive}
   ${\bf K^{\Delta}}$ is sound and strongly complete with respect to the class of all irreflexive frames. 
\end{theorem}

\begin{proof}
    By Thm.~\ref{thm.comp-k-delta}, it remains to show that
    \begin{center}
        $(\star)$~~~~~~every model has a pointwise-equivalent irreflexive model.
    \end{center}
    Given any model $\M=\lr{W,R,V}$, construct a model $\M'=\lr{W',R',V'}$ as follows.
\begin{itemize}
    \item[-] $W'=W\times \{+,-\}$, where we write $w^+$ for $(w,+)$ and $w^-$ for $(w,-)$,
    \item[-] for every $x,y\in S'$, $xR'y$ iff \\
    \begin{tabular}{ll}
     {\em either}    &  (i) $x=w^+$, $y=w^-$ and $wRw$\\
      {\em or}    &  (ii) $x=w^-$, $y=w^+$ and $wRw$\\
      {\em or} & (iii) $x=w^+$ or $x=w^-$, and $y=u^+$ or $y=u^-$, and $w\neq u$, and $wRu$.\\
    \end{tabular}
\item[-] for any $p\in\BP$ and any $x\in W'$, if either $x=w^+\text{ or }x=w^-$, then $V'(p)=\{x\mid w\in V(p)\}$.
\end{itemize}
It should be obvious that $\M'$ is irreflexive. In what follows, we show that for all $\phi\in\mathcal{L}(\Delta)$, for all $w\in W$, 
\begin{center}
($\heartsuit$)~~~~~~$\M,w\vDash\phi$ iff $\M',w^+\vDash\phi$ iff $\M',w^-\vDash\phi$.
\end{center}
We proceed by induction on $\phi$. The Boolean cases are straightforward. For the case $\Delta\phi$, suppose that $\M,w\nvDash\Delta\phi$, then there are $u,v$ such that $wRu$ and $wRv$ and $\M,u\vDash\phi$ and $\M,v\nvDash\phi$ (thus $u\neq v$). We consider three cases.
\begin{itemize}
    \item $w=u$, and thus $w\neq v$. Thus $\M,w\vDash\phi$. By induction hypothesis, $\M',w^+\vDash\phi$ and $\M',w^-\vDash\phi$. As $wRu$, we have $wRw$, which by definition of $R'$ implies that $w^+R'w^-$ and $w^-R'w^+$. As $wRv$ and $w\neq v$, then $w^+R'v^+$ (and $w^+R'v^-$) and $w^-R'v^+$ (and $w^-R'v^-$). By $\M,v\nvDash\phi$ and induction hypothesis, $\M',v^+\nvDash\phi$ (and $\M',v^-\nvDash\phi$). Therefore, $\M',w^+\nvDash\Delta\phi$ and $\M',w^-\nvDash\Delta\phi$.
    \item $w=v$, and thus $w\neq u$. Similar to the first case, we can show that $\M',w^+\nvDash\Delta\phi$ and $\M',w^-\nvDash\Delta\phi$.
    \item $w\neq u$ and $w\neq v$. Then by definition of $R'$, $w^+R'u^+$ and $w^+R'v^+$ and $w^-R'u^+$ and $w^-R'v^+$. By induction hypothesis, $\M',u^+\vDash\phi$ and $\M',v^+\nvDash\phi$. Therefore, $\M',w^+\nvDash\Delta\phi$ and $\M',w^-\nvDash\Delta\phi$.
\end{itemize}
Conversely, assume that $\M',w^+\nvDash \Delta\phi$ ($\M',w^-\nvDash\Delta\phi$, resp.). Then there are $y,z\in W'$ such that $w^+R'y$ ($w^-R'y$, resp.) and $w^+R'z$ ($w^-R'z$, resp.) and $\M',y\vDash\phi$ and $\M',z\nvDash\phi$ (thus $y\neq z$). We distinguish three cases.
\begin{itemize}
    \item $y=w^-$ ($y=w^+$, resp.). Then $z\neq w^-$ ($z\neq w^+$, resp.). Since $\M'$ is irreflexive, $z\neq w^-$ and $z\neq w^+$. As $w^+R'y$ ($w^-R'y$, resp.), we have $wRw$; as $w^+R'z$ ($w^-R'z$, resp.), according to the definition of $R'$, there must exist $u$ such that $w\neq u$ and $wRu$ and ($z=u^+$ or $z=u^-$). By induction hypothesis, $\M,u\nvDash\phi$. From $\M',y\vDash\phi$, $y=w^-$ ($y=w^+$, resp.) and induction hypothesis, it follows that $\M,w\vDash\phi$. Thus $\M,w\nvDash\Delta\phi$.
    \item $z=w^-$ ($z=w^+$, resp.). Similar to the first case, we can show that $\M,w\nvDash\Delta\phi$.
    \item $y\neq w^-$ and $z\neq w^-$ ($y\neq w^+$ and $z\neq w^+$, resp.). Since $\M'$ is irreflexive, we have also $y\neq w^+$ and $z\neq w^+$ ($y\neq w^-$ and $z\neq w^-$, resp.). Then according to the definition of $R'$, there must exist $u$ such that $wRu$ and $w\neq u$ and ($y=u^+$ or $y=u^-$), and $v$ such that $wRv$ and $w\neq v$ and ($z=v^+$ or $z=v^-$). By induction hypothesis, $\M,u\vDash\phi$ and $\M,v\nvDash\phi$. Therefore, $\M,w\nvDash\Delta\phi$.
\end{itemize}
This completes the proof of $(\heartsuit)$, thus we show that $(\star)$, as desired.
\end{proof}

Similar to the $(\star)$ in the proof of Thm.~\ref{thm.comp-kdelta-irreflexive}, we can also show that every symmetric model can be transformed into a
pointwise-equivalent
irreflexive symmetric model.


\section{Axiomatizations of $\mathcal{L}(\trianglecdot)$}\label{sec.axiomatizations-trianglecdot}

\begin{definition}[${\bf K^{\trianglecdot}}$ and ${\bf KB^{\trianglecdot}}$]\label{def.trianglecdot-systems} The proof system ${\bf K^{\trianglecdot}}$ consists of the following axioms and inference rules. 
\[
\begin{array}{ll}
   \text{TAUT}  & \text{all instances of propositional tautologies} \\
   \trianglecdot\text{Con} & \trianglecdot(\chi\to\phi)\land\trianglecdot(\neg\chi\to\phi)\to \trianglecdot\phi\\
   \trianglecdot\text{Dis} & \trianglecdot\phi\to \trianglecdot(\phi\to\psi)\vee\trianglecdot(\neg\phi\to\chi)\\
   \trianglecdot{\lra}  & \trianglecdot\phi\lra \trianglecdot\neg\phi\\
   \text{MP} & \text{From }\phi\text{ and }\phi\to\psi\text{ infer }\psi\\
   \text{NEC}\trianglecdot & \text{From }\phi\text{ infer }\trianglecdot\phi\\
   \text{RE}\trianglecdot & \text{From }\phi\lra \psi\text{ infer }\trianglecdot\phi\lra \trianglecdot\psi\\
\end{array}
\]
The proof system ${\bf KB^{\trianglecdot}}$ is the extension of ${\bf K^{\trianglecdot}}$ with axiom $\trianglecdot\text{B}$: $$\phi\to\trianglecdot((\trianglecdot\phi\land\trianglecdot(\phi\to\psi)\land\neg\trianglecdot\psi)\to\chi).$$
\end{definition}

We can see that Definitions~\ref{def.Delta-systems} and \ref{def.trianglecdot-systems} are the same but we only replace $\Delta$ by $\trianglecdot$ everywhere.

Notions of derivation and theorem are defined as usual. In particular, we use $\Gamma\vdash_{\bf K^{\trianglecdot}}\phi$ (resp. $\Gamma\vdash_{\bf KB^{\trianglecdot}}\phi$) to denote that $\phi$ is derivable from $\Gamma$ in ${\bf K^{\trianglecdot}}$ (resp. ${\bf KB^{\trianglecdot}}$).

In what follows, we reduce the soundness and strong completeness of ${\bf K^{\trianglecdot}}$ to those of ${\bf K^{\Delta}}$, via a pair of translation functions, namely $\trianglecdot$-translation $(\cdot)^{\trianglecdot}$ and $\Delta$-translation $(\cdot)^{\Delta}$ (defined below). The strategy can be summarized as follows. For all $\Gamma\cup\{\phi\}\subseteq\mathcal{L}(\trianglecdot)$,
$$\Gamma\vdash_{\bf K^{\trianglecdot}}\phi\stackrel{(1)}\iff \Gamma^{\Delta}\vdash_{\bf K^{\Delta}}\phi^{\Delta}\stackrel{(2)}\iff \Gamma^{\Delta}\vDash_{K}\phi^{\Delta}\stackrel{(3)}\iff \Gamma\vDash_{K}\phi,$$
$$\Gamma\vdash_{\bf KB^{\trianglecdot}}\phi\iff \Gamma^{\Delta}\vdash_{\bf KB^{\Delta}}\phi^{\Delta}\iff \Gamma^{\Delta}\vDash_{B}\phi^{\Delta}\iff \Gamma\vDash_{B}\phi.$$

In what follows, we will show the case for ${\bf K^{\trianglecdot}}$ in detail. The proof for ${\bf KB^{\trianglecdot}}$ is similar.

\medskip

We now define two translations, namely, $\trianglecdot$-translation and $\Delta$-translation. The {\em $\trianglecdot$-translation}, denoted $(\cdot)^{\trianglecdot}$, replaces every occurrence of $\Delta$ in every $\mathcal{L}(\Delta)$-formula with $\trianglecdot$; the {\em $\Delta$-translation}, denoted $(\cdot)^{\Delta}$, replaces every occurrence of $\trianglecdot$ in every $\mathcal{L}(\trianglecdot)$-formula with $\Delta$. It is straightforward to show by induction that $(\phi^{\trianglecdot})^{\Delta}=\phi$ for all $\phi\in\mathcal{L}(\Delta)$, and $(\psi^{\Delta})^{\trianglecdot}=\psi$ for all $\psi\in\mathcal{L}(\trianglecdot)$. In general, $(\Gamma^{\trianglecdot})^{\Delta}=\Gamma$ for all $\Gamma\subseteq\mathcal{L}(\Delta)$ and $(\Sigma^{\Delta})^{\trianglecdot}=\Sigma$ for all $\Sigma\subseteq\mathcal{L}(\trianglecdot)$.

Due to the purely notational difference between the proof systems ${\bf K^{\trianglecdot}}$ and ${\bf K^{\Delta}}$, we have the following.
\begin{lemma}\label{lem.notationaldifference}
For all $\Gamma\cup\{\phi\}\subseteq\mathcal{L}(\trianglecdot)$,
$$\Gamma\vdash_{\bf K^{\trianglecdot}}\phi\iff\Gamma^{\Delta}\vdash_{\bf K^{\Delta}}\phi^{\Delta}.$$
\end{lemma}

This completes the step (1) in the above strategy. The step (2) follows directly from 
Thm.~\ref{thm.comp-k-delta}. So it remains only to show the step (3).  For this, we need some preparations.

Recall that we have remarked that the semantics of $\trianglecdot$ and $\Delta$ are the same on the irreflexive models. Here is a precise statement of that result, which can be proved by induction on $\mathcal{L}(\trianglecdot)$-formulas.
\begin{proposition}\label{prop.irrelfexive-equivalent}
    For all irreflexive models $\M$, for all worlds $w$ in $\M$, for all $\phi\in\mathcal{L}(\trianglecdot)$, we have
    $$\M,w\vDash \phi\iff \M,w\vDash \phi^{\Delta}.$$
Consequently, for all $\Gamma\subseteq\mathcal{L}(\trianglecdot)$, $\M,w\vDash\Gamma\iff \M,w\vDash\Gamma^{\Delta}$.
\end{proposition}

Recall from Def.~\ref{def.ir} that $\M^{-T}$ is the irreflexive reduction of $\M$, thus $\M^{-T}$ is irreflexive. Then the result below follows directly from Prop.~\ref{prop.invariant-ir} and Prop.~\ref{prop.irrelfexive-equivalent}.
\begin{corollary}\label{coro.iff}
    For all models $\M$, for all $w\in W$, for all $\phi\in\mathcal{L}(\trianglecdot)$, we have
    $$\M,w\vDash\phi\iff \M^{-T},w\vDash\phi^{\Delta}.$$
Consequently, for all $\Gamma\subseteq\mathcal{L}(\trianglecdot)$, $\M,w\vDash\Gamma\iff \M^{-T},w\vDash\Gamma^{\Delta}$.
\end{corollary}

\weg{\begin{proposition}
    ${\bf K^{\trianglecdot}}$ is sound with respect to the class of all frames.
\end{proposition}

Recall from ... that $\M^{-T}$ is the irreflexive reduction of $\M$.}

\begin{lemma}\label{lem.equivalence}
    For all $\Gamma\cup\{\phi\}\subseteq\mathcal{L}(\trianglecdot)$, we have
    $$\Gamma\vDash_{K}\phi\iff \Gamma^{\Delta}\vDash_{K}\phi^{\Delta}.$$
\end{lemma}

\begin{proof}
    Suppose that $\Gamma\nvDash_{K}\phi$. Then there exists a model $\M$ and a world $w$ in $\M$ such that $\M,w\vDash\Gamma$ and $\M,w\nvDash\phi$. By Coro.~\ref{coro.iff}, $\M^{-T},w\vDash\Gamma^{\Delta}$ and $\M^{-T},w\nvDash\phi^{\Delta}$. Thus $\Gamma^{\Delta}\nvDash_{K}\phi^{\Delta}$.

    Conversely, assume that $\Gamma^{\Delta}\nvDash_{K}\phi^{\Delta}$. Then there exists a model $\M$ and a world $w$ in $\M$ such that $\M,w\vDash\Gamma^{\Delta}$ and $\M,w\nvDash\phi^{\Delta}$. Due to $(\ast)$ in the proof of Thm.~\ref{thm.comp-kdelta-irreflexive}, there exists an irreflexive model $\M'$ and $w'$ in $\M'$ such that $\M',w'\vDash\Gamma^{\Delta}$ and $\M',w'\nvDash\phi^{\Delta}$. Now using Prop.~\ref{prop.irrelfexive-equivalent}, we infer that $\M',w'\vDash\Gamma$ and $\M',w'\nvDash\phi$. Therefore, $\Gamma\nvDash_{K}\phi$.
\end{proof}

This finishes the step (3) in the above strategy. Based on the previous analysis, we obtain the soundness and strong completeness of ${\bf K^{\trianglecdot}}$. As a matter of fact, we can establish the following general result.
\begin{theorem}
    Let $T\subseteq X$. Then ${\bf K^{\trianglecdot}}$ is sound and strongly complete with respect to $X$. That is, for all $\Gamma\cup\{\phi\}\subseteq\mathcal{L}(\trianglecdot)$, we have:
    $$\Gamma\vdash_{\bf K^{\trianglecdot}}\phi\iff \Gamma\vDash_{X}\phi.$$
\end{theorem}

\begin{proof}
    Let $\Gamma\cup\{\phi\}\subseteq\mathcal{L}(\trianglecdot)$. We have the following equivalences.
    $$\Gamma\vdash_{\bf K^{\trianglecdot}}\phi\iff \Gamma^{\Delta}\vdash_{\bf K^{\Delta}}\phi^{\Delta}\iff\Gamma^{\Delta}\vDash_{K}\phi^{\Delta}\iff\Gamma\vDash_{K}\phi\iff\Gamma\vDash_{X}\phi,$$
    where the first equivalence follows from Lemma~\ref{lem.notationaldifference}, the second from Thm.~\ref{thm.comp-k-delta}, the third from Lemma~\ref{lem.equivalence}. For the last equivalence, note that $T\subseteq K$ and $T\subseteq X$. By Prop.~\ref{prop.semanticproperties}(1), $\Gamma\vDash_{K}\phi$ iff $\Gamma\vDash_{T}\phi$, and $\Gamma\vDash_{X}\phi$ iff $\Gamma\vDash_{T}\phi$, and thus $\Gamma\vDash_{K}\phi$ iff $\Gamma\vDash_{X}\phi$.
\end{proof}

\begin{corollary}\label{coro.k-trianglecdot}
${\bf K^{\trianglecdot}}$ is sound and strongly complete with respect to the class $K$ of all frames, to the class $D$ of serial frames, and to the class $T$ of reflexive frames.     
\end{corollary}

As with the reduction of the sound and strong completeness of ${\bf K^{\trianglecdot}}$ to those of
${\bf K^{\Delta}}$, we can also reduce the soundness and strong completeness of ${\bf KB^{\trianglecdot}}$ to those of ${\bf KB^{\Delta}}$.
Note that every symmetric model can be transformed into a
pointwise-equivalent
irreflexive symmetric model (see the remark after Thm.~\ref{thm.comp-kdelta-irreflexive}), thus we can show a similar result
to Lemma~\ref{lem.equivalence}. Then, according to the previous analysis, we can obtain the soundness and strong completeness of ${\bf KB^{\trianglecdot}}$ from those of ${\bf KB^{\Delta}}$ given in~\cite{Fanetal:2014}. Once again, as a matter of fact, we have the following general result.
\begin{theorem}
    Let ${TB}\subseteq X\subseteq B$. Then ${\bf KB^{\trianglecdot}}$ is sound and strongly complete with respect to $X$. That is, for all $\Gamma\cup\{\phi\}\subseteq\mathcal{L}(\trianglecdot)$, we have:
    $$\Gamma\vdash_{\bf KB^{\trianglecdot}}\phi\iff \Gamma\vDash_{X}\phi.$$
\end{theorem}

\begin{corollary}\label{coro.kb-trianglecdot}
    ${\bf KB^{\trianglecdot}}$ is sound and strongly complete with respect to the class $B$ of symmetric frames, to the class $DB$ of serial and symmetric frames, and also to the class $TB$ of reflexive and symmetric frames.
\end{corollary}

We have seen that all axioms of ${\bf K^{\trianglecdot}}$ and axiom $\trianglecdot\text{B}$ are similar in form to their counterparts in $\mathcal{L}(\Delta)$, respectively. It would be natural to ask whether we can obtain axioms over other special frame classes by translating corresponding axioms in $\mathcal{L}(\Delta)$ via the $\trianglecdot$-translation. As we shall see, the answer would be negative. 
\begin{proposition}\
 \begin{itemize}
     \item[(1)] Over the class of reflexive frames, $\Delta\phi\land\Delta(\phi\to\psi)\land\phi\to\Delta\psi$ is valid, but $\trianglecdot\phi\land\trianglecdot(\phi\to\psi)\land\phi\to\trianglecdot\psi$ is not.
     \item[(2)] Over the class of transitive frames, $\Delta\phi\to \Delta(\Delta\phi\vee\psi)$ is valid, but $\trianglecdot\phi\to \trianglecdot(\trianglecdot\phi\vee\psi)$ is not.
     \item[(3)] Over the class of Euclidean frames, $\neg\Delta\phi\to \Delta(\neg\Delta\phi\vee\psi)$ is valid, but $\neg\trianglecdot\phi\to \trianglecdot(\neg\trianglecdot\phi\vee\psi)$ is not. 
 \end{itemize}   
\end{proposition}

\begin{proof}
  The frame validities can be found in~\cite[Prop.~5.2]{Fanetal:2015}. We only need to consider the frame invalidities.
  \begin{itemize}
     \item[(1)] Consider the following reflexive model $\M_1$ (where other valuations are inessential):
     
     \[
     \xymatrix{u_1:\neg p,q\ar@(ul,ur)&w_1:p\ar[l]\ar[r]\ar@(ul,ur)&v_1:\neg p,\neg q\ar@(ur,ul)\\}
     \]
     Then $\M_1,w_1\vDash \trianglecdot p$ and $\M_1,w_1\vDash\trianglecdot(p\to q)$ and $\M_1,w_1\vDash p$ but $\M_1,w_1\nvDash\trianglecdot q$. Therefore, $\M_1,w_1\nvDash\trianglecdot p\land\trianglecdot(p\to q)\land p\to\trianglecdot q$.
     \item[(2)] Consider the following transitive model $\M_2$ (where other valuations are inessential):

     \[
     \xymatrix{&u_2:\neg p\\
     w_2:p\ar[r]\ar[ur]& v_2:\neg p,\neg q\ar[l]\ar[u]\\}
     \]
    Then that $\M_2,w_2\vDash \trianglecdot p$ and $\M_2,w_2\nvDash\trianglecdot(\trianglecdot p\vee q)$ (since $u_2\vDash\trianglecdot p$, we have $\M_2,u_2\vDash \trianglecdot p\vee q$; however, $\M_2,v_2\nvDash\trianglecdot p\vee q$). Therefore, $\M_2,w_2\nvDash\trianglecdot p\to \trianglecdot(\trianglecdot p\vee q)$.
     \item[(3)] Consider the following Euclidean model $\M_3$ (where other valuations are inessential):

     \[
     \xymatrix{
     & u_3: p,\neg q\ar[d]\ar@(ur,ul)\\
     w_3\ar[ur]\ar[r]&v_3:\neg p,q\ar[u]\ar@(dr,dl)\\
     }
     \]
     
     Then $\M_3,w_3\vDash\neg\trianglecdot p$ and $\M_3,w_3\nvDash\trianglecdot(\neg\trianglecdot p\vee q)$. To see the latter, note that $\M_3,u_3\nvDash\neg\trianglecdot p\vee q$ and $\M_3,v_3\vDash \neg\trianglecdot p\vee q$. Therefore, $\M_3,w_3\nvDash \neg\trianglecdot p\to \trianglecdot(\neg\trianglecdot p\vee q)$.
 \end{itemize}
\end{proof}

\begin{corollary}\label{coro.necessary}
    Let $\mathfrak{C}$ be the class of reflexive frames, the class of transitive frames, or the class of Euclidean frames. Then there is a $\phi\in\mathcal{L}(\Delta)$ such that $\mathfrak{C}\vDash\phi$ but $\mathfrak{C}\nvDash\phi^{\trianglecdot}$.   
\end{corollary}

Thus all axioms of ${\bf K^{\trianglecdot}}$ are valid over the class of all frames, and axiom $\trianglecdot\text{B}$ is valid over the class of symmetric frames; however, it is not the case for other special axioms. One may wonder why this is the case. In what follows, we explain this phenomenon.
\begin{proposition}\label{prop.sufficient}
    Let $\mathfrak{C}$ be a class of frames. If (a) $\mathfrak{C}$ is closed under irreflexivization, then (b) for all $\phi\in\mathcal{L}(\Delta)$, if $\mathfrak{C}\vDash\phi$, then $\mathfrak{C}\vDash\phi^{\trianglecdot}$.
\end{proposition}

\begin{proof}
    Suppose that (a) holds, to show that (b) holds. For this, let $\phi\in\mathcal{L}(\Delta)$ and assume that $\mathfrak{C}\nvDash\phi^{\trianglecdot}$, it remains only to prove that $\mathfrak{C}\nvDash\phi$.

    By assumption, there is a $\mathfrak{C}$-model $\M$ and a world $w$ in $\M$ such that $\M,w\nvDash\phi^{\trianglecdot}$. Since $\phi^{\trianglecdot}\in\mathcal{L}(\trianglecdot)$, it follows from Coro.~\ref{coro.iff} that $\M^{-T},w\nvDash(\phi^{\trianglecdot})^{\Delta}$, that is, $\M^{-T},w\nvDash\phi$. As $\M$ is a $\mathfrak{C}$-model, by (a), $\M^{-T}$ is also a $\mathfrak{C}$-model. Therefore, $\mathfrak{C}\nvDash\phi$.
\end{proof}

Since all axioms of ${\bf K^{\Delta}}$ and axiom $\Delta\text{B}$ are valid over the class of all frames and over the class of symmetric frames~\cite{Fanetal:2015}, respectively, and both frame classes are closed under irreflexivization, by Prop.~\ref{prop.sufficient}, all axioms of ${\bf K^{\trianglecdot}}$ and axiom $\trianglecdot\text{B}$ are valid over the corresponding classes of frames.
\weg{\begin{corollary}\
    \begin{itemize}
    \item[(1)] All axioms of ${\bf K^{\trianglecdot}}$ are valid over the class of all frames.
    \item[(2)] $\phi\to\trianglecdot(\trianglecdot\phi\land\trianglecdot(\phi\to\psi)\land\neg\trianglecdot\psi\to\chi)$ is valid over the class of symmetric frames.
    \end{itemize}
\end{corollary}}

It follows from Prop.~\ref{prop.sufficient} that (a) is a {\em sufficient} condition of (b). A natural question may be whether (a) is also a {\em necessary} condition of (b). The answer is positive when $\mathfrak{C}$ is the class of reflexive frames, the class of transitive frames, or the class of Euclidean frames, which follows from Coro.~\ref{coro.necessary} and the fact that none of such classes of frames is closed under irreflexivization (for the case of the class of reflexive frames, consider a single-point reflexive frame; for the case of the class of transitive frames or the class of Euclidean frames, consider a two-point universal frame). However, in general, the answer would be negative.
\begin{proposition}
    Let $\mathfrak{C}$ be the class of all frames where each frame has at least one reflexive point. Then (a) of Prop.~\ref{prop.sufficient} fails, but (b) of Prop.~\ref{prop.sufficient} holds.
\end{proposition}

\begin{proof}
    One may check that $\mathfrak{C}$ is not closed under irreflexivization, thus (a) of Prop.~\ref{prop.sufficient} fails. To show (b) of Prop.~\ref{prop.sufficient}, we will show that for all $\phi\in\mathcal{L}(\Delta)$, (1) $\mathfrak{C}\vDash\phi$ implies $\vDash\phi$, (2) $\vDash\phi$ implies $\vDash\phi^{\trianglecdot}$, and (3) $\vDash\phi^{\trianglecdot}$ implies $\mathfrak{C}\vDash\phi^{\trianglecdot}$. As (3) is obvious, it remains only to show (1) and (2).

    For (1), suppose that $\nvDash\phi$, then $\neg\phi$ is satisfiable. 
    Then there exists a model $\M$ and a world $w$ in $\M$ such that $\M,w\vDash\neg\phi$. Let $w'\notin \M$.\footnote{If $\M$ is finite, then such a $w'$ must exist. Even if $\M$ is infinite, we can also add more members --- remember the Hilbert Hotel!} Define $\M'=\lr{W',R',V'}$, where $W'=\{w'\}$ and $R'=\{(w',w')\}$. Now consider the disjoint union of $\M$ and $\M'$, which is denoted as $\M''$. Since $\M''$ contains a reflexive point, its underlying frame belongs to $\mathfrak{C}$. Also, as $\mathcal{L}(\Delta)$ is a fragment of $\mathcal{L}(\Box)$ and the satisfaction of $\mathcal{L}(\Box)$-formulas is invariant under disjoint union (e.g.~\cite{blackburnetal:2001}), so is the satisfaction of $\mathcal{L}(\Delta)$-formulas. This follows that $\M'',w\vDash\neg\phi$, and thus $\mathfrak{C}\nvDash\phi$.

    For (2), assume that $\nvDash\phi^{\trianglecdot}$, then there exists a model $\mathcal{N}$ and a world $w$ in $\mathcal{N}$ such that $\mathcal{N},w\nvDash\phi^{\trianglecdot}$. Since $\phi^{\trianglecdot}\in \mathcal{L}(\trianglecdot)$, by Coro.~\ref{coro.iff} we infer that $\mathcal{N}^{-T},w\nvDash(\phi^{\trianglecdot})^{\Delta}$, that is, $\N^{-T},w\nvDash\phi$. Therefore, $\nvDash\phi$.
\end{proof}

\section{Conclusion and Future Work}\label{sec.conclusion}

In this paper, inspired by the notions of contingency and strong possibility, we proposed a notion called `strong contingency', which can be seen as a mixture of the former two notions. We proposed the bisimulation for the logic $\mathcal{L}(\trianglecdot)$ with this notion as a sole primitive modality. We investigated the expressivity and frame definability of $\mathcal{L}(\trianglecdot)$. It turns out that $\mathcal{L}(\trianglecdot)$ and the logic of contingency $\mathcal{L}(\Delta)$ (and also the standard modal logic $\mathcal{L}(\Box)$) are incomparable in expressivity even over the class of $\mathcal{S}5$-models, and $\mathcal{L}(\trianglecdot)$ is less expressive than $\mathcal{L}(\boxdot)$ over the class of $\mathcal{S}5$-models. Last but not least, we completely axiomatized $\mathcal{L}(\trianglecdot)$ over various frame classes, where the soundness and strong completeness of the axiomatizations are reduced to those of the corresponding axiomatizations of $\mathcal{L}(\Delta)$.

Several issues remain for future work. For instance, one may axiomatize $\mathcal{L}(\trianglecdot)$ over other special frame classes, such as the class of transitive frames and the class of Euclidean frames. As we demonstrated in Prop.~\ref{prop.almost-definability}, $\boxdot$ is almost definable in terms of $\trianglecdot$, which indicates that we may use the method based on almost-definability schemas proposed in~\cite{Fan:2025bundled}. For another example, one may investigate the bimodal logic with $\trianglecdot$ and $\Delta$, and the bimodal logic with $\trianglecdot$ and $\Box$. Moreover, one may apply the bifurcation method to explore the variants of the operator $\trianglecdot$ and axiomatize them over various frame classes, as with the case of the non-contingency operator in~\cite{Fan:2021twovariants}.

\paragraph*{Acknowledgements}

This article is supported by the Key Project of National Social Science Fund of China (no. 24AZX019). The author is grateful to Lloyd Humberstone for drawing this topic to his attention. The author thanks an anonymous reviewer for his/her insightful comments.

\bibliographystyle{plain}
\bibliography{sample,biblio2024}

\begin{thebibliography}{10}

\bibitem{blackburnetal:2001}
P.~Blackburn, M.~de~Rijke, and Y.~Venema.
\newblock {\em Modal Logic}.
\newblock Cambridge University Press, Cambridge, 2001.
\newblock Cambridge Tracts in Theoretical Computer Science 53.

\bibitem{Borgan67}
A.~Brogan.
\newblock Aristotle's logic of statements about contingency.
\newblock {\em Mind}, 76(301):49--61, 1967.

\bibitem{Chen:2020}
J.~Chen.
\newblock Logics for strong possibility and weak necessity.
\newblock {\em Studies in Logic}, 13(1):53--74, 2020.
\newblock {(In Chinese)}.

\bibitem{Cresswell88}
M.~Cresswell.
\newblock Necessity and contingency.
\newblock {\em Studia Logica}, 47:145--149, 1988.

\bibitem{deRijke:1992}
M.~de~Rijke.
\newblock The modal logic of inequality.
\newblock {\em Journal of Symbol Logic}, 57:566--584, 1992.

\bibitem{Fan:2019symmetric}
J.~Fan.
\newblock Symmetric contingency logic with unlimitedly many modalities.
\newblock {\em Journal of Philosophical Logic}, 48(5):851–866, 2019.

\bibitem{Fan:2021twovariants}
J.~Fan.
\newblock Two variants of noncontingency operator.
\newblock {\em Notre Dame Journal of Formal Logic}, 62(3):491--525, 2021.

\bibitem{Fan:2023}
J.~Fan.
\newblock Logics of non-actual possible worlds.
\newblock {\em Studies in Logic}, 16(6):17–38, 2023.

\bibitem{Fan:2025bundled}
J.~Fan.
\newblock A uniform approach to axiomatizing bundled operators.
\newblock {\em The Review of Symbolic Logic}, 18(4):1121--1167, 2025.
\newblock https://doi.org/10.1017/S1755020325100981.

\bibitem{Fanetal:2014}
J.~Fan, Y.~Wang, and H.~van Ditmarsch.
\newblock Almost necessary.
\newblock In {\em Advances in Modal Logic}, volume~10, pages 178--196, 2014.

\bibitem{Fanetal:2015}
J.~Fan, Y.~Wang, and H.~van Ditmarsch.
\newblock Contingency and knowing whether.
\newblock {\em The Review of Symbolic Logic}, 8(1):75--107, 2015.

\bibitem{Humberstone95}
L.~Humberstone.
\newblock The logic of non-contingency.
\newblock {\em Notre Dame Journal of Formal Logic}, 36(2):214--229, 1995.

\bibitem{Humberstone:2007}
L.~Humberstone.
\newblock Modal logic for other-world agnostics: neutrality and halld\'{e}n
  incompleteness.
\newblock {\em Journal of Philosophical Logic}, 36(1):1--32, 2007.

\bibitem{DBLP:journals/ndjfl/Kuhn95}
S.~Kuhn.
\newblock Minimal non-contingency logic.
\newblock {\em Notre Dame Journal of Formal Logic}, 36(2):230--234, 1995.

\bibitem{Malink:2016}
M.~Malink.
\newblock Aristotle on one-sided possibility.
\newblock In Rini M.~Cresswell, A and E.~Mares, editors, {\em Logical
  Modalities from Aristotle to Carnap: The Story of Necessity}, pages 29--49.
  Cambridge University Press, 2016.

\bibitem{MR66}
H.~Montgomery and R.~Routley.
\newblock Contingency and non-contingency bases for normal modal logics.
\newblock {\em Logique et Analyse}, 9:318--328, 1966.

\bibitem{MR68}
H.~Montgomery and R.~Routley.
\newblock Non-contingency axioms for {S4 and S5}.
\newblock {\em Logique et Analyse}, 11(43):422--424, 1968.

\bibitem{Segerberg:1976}
K.~Segerberg.
\newblock {``Somewhere else'' and ``Some other time''}.
\newblock In {\em Wright and Wrong: Mini-Essays in Honor of George Henrik von
  Wright on his Sixtieth Birthday}, pages 61--64. Publications of the Group in
  Logic and Methodology of Real Finland, Vol. 3, {\AA}bo Akademi, 1976.

\bibitem{Segerberg:1980}
K.~Segerberg.
\newblock A note on the logic of elsewhere.
\newblock {\em Theoria}, 46:183--187, 1980.

\bibitem{hoeketal:2004}
W.~van~der Hoek and A.~Lomuscio.
\newblock A logic for ignorance.
\newblock {\em Electronic Notes in Theoretical Computer Science},
  85(2):117--133, 2004.

\bibitem{Wright:1979}
G.~H. von Wright.
\newblock A modal logic of place.
\newblock In E.~Sosa, editor, {\em The Philosophy of Nicholas Rescher:
  Discussions and Replies}, pages 65--73. Reidel, Dordrecht, 1979.

\bibitem{DBLP:journals/ndjfl/Zolin99}
E.~Zolin.
\newblock Completeness and definability in the logic of noncontingency.
\newblock {\em Notre Dame Journal of Formal Logic}, 40(4):533--547, 1999.

\end{thebibliography}

\end{document}